\documentclass[10pt,reqno]{amsart}
\usepackage{latexsym}
\usepackage{graphicx}
\graphicspath{ {./images/} }
\usepackage{fancyhdr,amssymb}
\usepackage{epsfig}

\theoremstyle{plain}
\numberwithin{equation}{section}
\newtheorem{thm}{Theorem}[section]

\theoremstyle{definition}

\begin{document}
\fancyhead{}
\renewcommand{\headrulewidth}{0pt}
\fancyfoot{}
\fancyfoot[LE,RO]{\medskip \thepage}

\setcounter{page}{1}

\title[Sunlet Factors for Cartesian Products of Cycles]{Sunlet Factors for Cartesian Products of Cycles}
\author{Henry Jervis}
\address{Dept of Mathematics and Statistics\\
                Georgetown University\\
                Washington D.C.\\
                USA}
\email{hwj7@georgetown.edu}
\thanks{}
\author{Paul C. Kainen}
\address{Dept of Mathematics and Statistics\\
                Georgetown University\\
                Washington D.C.\\
                USA}
\email{kainen@georgetown.edu}

\begin{abstract}
A sunlet is a cycle with a pendant edge attached at each vertex of the cycle. For the bipartite toroidal grid graphs $C_{2n} \,\Box\,C_{2n}$, factorizations into sunlets are given by homomorphisms from disjoint unions of $s$ copies of a sunlet for $s \in \{1, n, n^2\}$, $n \geq 3$, such that edges are mapped bijectively.
\end{abstract}

\maketitle

\section{Introduction}
The problem of factoring a graph $G$ into isomorphic edge-disjoint copies of some smaller graph $H$ has been considered (e.g., Bos\'ak  \cite{bosak}) for a variety of pairs $(G,H)$.  For example, if $G=K_n$ and $H= C_n$, then the problem is solvable for $n \geq 3$ odd; see Alspach \cite{alspach}.  This corresponds to Hamiltonian cycle factorizations of the complete graph and is a classical problem in recreational mathematics, involving seatings around a table with no repeated neighbors. 

Another example has $G=K_n$ and $H=Q_5$, the 5-dimensional cube, where Bryant et al \cite{bemv} showed the necessary divisibility conditions on $n$, that 5 divides $n$ and that 80 divides $n(n-1)/2$, are also sufficient. 

Quite generally, El-Zanati and Vanden Eynden \cite{ev} enumerated the factorizations of cycle products by cycles of all lengths; those allowed by divisibility all exist.

 Rather than requiring  $H$ to be regular (all vertices of equal degree), Anitha and Lekshmi \cite{anitha} considered decompositions into {\it $n$-sunlets} (they called them ''suns'') which are $n$-cycles, with an additional set of $n$ endpoint vertices such that each vertex in the cycle is adjacent to a unique endpoint. \par
\ \\
\noindent\rule{0.84in}{0.4pt} \par
\medskip
\indent\indent {\fontsize{8pt}{9pt} \selectfont MSC2020: 05C60, 05C38, 05C62 \par}
\indent\indent {\fontsize{8pt}{9pt} \selectfont Key words and phrases: cycle, covering, factorization, finite state machine. \par}

\thispagestyle{fancy}

\vfil\eject
\fancyhead{}
\fancyhead[CO]{\hfill SUNLET FACTORS FOR CARTESIAN PRODUCTS OF CYCLES}
\fancyhead[CE]{H.~JERVIS AND P.~KAINEN  \hfill}
\renewcommand{\headrulewidth}{0pt}

An $n$-sunlet has $2n$ vertices and $2n$ edges, where half the vertices have degree 3 and the other half are endpoints; there is a unique $n$-cycle through the cubic nodes. The Harary square of $C_{2n}$ is factorable into $n$-sunlets  \cite{anitha}. Sowndhariya and Muthusamy \cite{sm2021} factored $G$, a Cartesian product of complete graphs, into copies of the sunlet of order 8. See \cite{muth} for $G$ a hypercube.  

We factor the Cartesian product $C_{2n} \,\Box \, C_{2n}$ (toroidal grids) into sunlets via homomorphisms from disjoint unions of $s$ isomorphic sunlets onto toroidal grids, for $s=1$, $s=n$, and $s=n^2$, $n \geq 3$ such that the sunlet cycles are mapped isomorphically.  Further, it is possible to orient both graphs so that the sunlets are FSMs (out-degree 1 at every vertex) and the embedding respects the digraphs.

An outline of the paper follows.  Section 2 has definitions; Section 3 has Theorem 1, where $s=1$; Section 4 has Theorem 2 dealing with $s = n$; in Section 5, we have Theorem 3 for $s = n^2$; Section 6 is a brief discussion.

\section{Definitions}

For graph definitions, see \cite{harary, diestel}.  
Let $C_p$ denote a cycle of length $p$.  An edge $e$ in a graph $G$ is \textbf{pendant} if exactly one endpoint of $e$ has degree 1.  For $p \geq 3$, the \textbf {sunlet} $S^1_p$ is a supergraph of $C_p$ which adds $p$ pendant edges; $S^1_p$ is denoted $L_{2p}$ in \cite{muth, sm2021}. We write \textbf {Z(S)} for the unique cycle contained in a sunlet $S$.

A \textbf {factorization} of graph $G$ is a set of pairwise isomorphic subgraphs such that each edge of $G$ is in a unique subgraph. This departs from the usage in \cite{harary}. The \textbf {size} (number of edges) of $S^1_p$ is $2p$.  As $S^1_p$ is unicyclic and connected, size equals \textbf {order} (number of vertices). Write $G \,\Box \, H$ for \textbf {Cartesian product} of graphs $G$ and $H$. Let $n \cdot G \;:= \; \overline{K}_n \,\Box \, G$ denote $n$ disjoint copies of $G$.

A Cartesian product of cycles is called a \textbf {toroidal grid} as it can be embedded in the torus such that all regions have 4 sides.  An \textbf {orientation} of a graph $G=(V,E)$ is a digraph $D=(V,A)$, where for each $uw \in E$ one has either $(u,w)$ or $(w,u)$ in $A$ (not both); i.e., each edge is assigned a unique \emph{direction}.  A toroidal grid has 4 isomorphic \textbf {canonical orientations}; for the annular picture of a torus, we take \emph{outward} for $y$ and \emph{clockwise} for $x$. See Fig. \ref{fig:coords}.

A \textbf {homomorphism} $\varphi: G \to H$ from $G$ to $H$ is a function 

\begin{equation}\label{1}
\varphi:V_G \to V_H
\end{equation}

\noindent such that if $e = uw \in E_G$, then $\varphi(u)\,\varphi(w) \in E_H$. Define $\varphi(e):= \varphi(u)\,\varphi(w)$; $\varphi$ is \textbf {onto} if for every $e' = u'w' \in E_H$, there is $e=uw \in E_G$ such that $e' = \varphi(e)$. A 
homomorphism $\varphi: G \to H$ is \textbf {r-to-1 (on 
vertices)} if $|\varphi^{-1}(u')| = r$ for every $u' \in V_H$, where $\varphi^{-1}(u') = \{u \in V_G: \varphi(u) = u'\}$.  

Given a homomorphism $\varphi: G \to H$ and orientations of $G$ and $H$, we say that $\varphi$ \textbf {is compatible with the orientations} if $(\varphi(u), \varphi(w))$ is an arc when $(u,w)$ is an arc.

An \textbf {epimorphism} is an onto homomorphism and a \textbf {covering} is an epimorphism which is \textbf {1-to-1 (on edges)} (i.e., the edge-set is mapped bijectively).  The coverings defined below provide ``rearrangements'' of the arcs of a directed cycle product as a factorization into sunlets that preserves orientation. 

An oriented graph is a Finite State Machine (\textbf {FSM}) if each vertex has a single arc from it.  A sunlet has an FSM orientation by making the cycle a directed cycle and having each ray oriented toward the cycle. See Fig. \ref{fig:s=1}.

\section{Homomorphism from One Sunlet}

The following theorem is slightly stronger than the abstract as the toroidal grid need not be bipartite; additional information is given on orientations.

\begin{thm}
For $n \geq 3$, there is a 2-to-1 covering
$$\varphi: S := S_{n^2}^{1} \to C_{n} \ \square \ C_{n}$$ 
such that $\varphi(Z(S))$ is a Hamiltonian cycle.  Further, $\varphi$ is compatible with the FSM orientation of $S$ and the standard orientation for the grid.
\end{thm}

\begin{proof}
    
Let each vertex in the graph $C_{n} \ \square \ C_{n}$ be denoted by $(x,y)$, where \textbf {x indexes the concentric n-gons} and \textbf {y is the location within the n-gon}, both from $[0,n-1]$, coordinates are mod $n$; see Fig. \ref{fig:coords}. Call the edges of $C_{n} \ \square \ C_{n}$ {\bf vertical} if they connect concentric n-gons and \textbf {horizontal} if they make up the n-gon itself. Embed the Hamiltonian cycle as shown in Fig. \ref{fig:s=1}; the figure also illustrates the rays incident to each cycle vertex.  This cycle follows $n{-}1$ consecutive horizontal edges before turning onto a vertical edge.

\begin{figure}[ht]
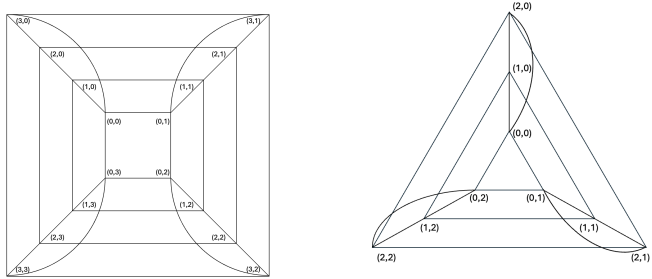

    \centering
    \includegraphics[scale = .2]{images/C4C4GraphwithCoordinates.pdf}
    \includegraphics[scale=.27]{images/C3C3GraphwithCoordinates.pdf}
    \caption{$C_{4} \ \square \ C_{4}\;$ and $\;C_{3} \ \square \ C_{3}\;$ with described coordinates.}
    \label{fig:coords}
\end{figure}

\noindent Hence, $\varphi: V(C_{n^2}) \to V(C_{n} \ \square \ C_{n})$ so that vertices $0, 1, \ldots, n^2 {-} 1, 0$ map to
\[
H = ((0,0),\ldots,(0,n{-}1),(1,n{-}1),(1,0),\ldots,(1,n{-}2),\ldots, (n{-}1,0),(0,0)).
\]

Since cycle $H$ spans the 4-regular graph $C_{n} \ \square \ C_{n}$, the remaining edges induce a degree-2 graph $H'$.  In each connected component of $H'$ orient the edges into a dicycle.  Each vertex $u$ of $C_{n} \ \square \ C_{n}$ is the target of a unique arrow $(v,u)$ in $H'$ and the ray $w'w$ at the cycle vertex $w \in \varphi^{-1}(u)$ of $S_{n^2}^1$ is mapped so that $\varphi(w',w) = (v,u)$.  This shows the existence of $\varphi$.  However, for algorithmic uses, we explicitly construct an extension where $C_{n} \ \square \ C_{n}$ has standard orientation.

\begin{figure}[ht]
    \centering
    \includegraphics[scale=.3]{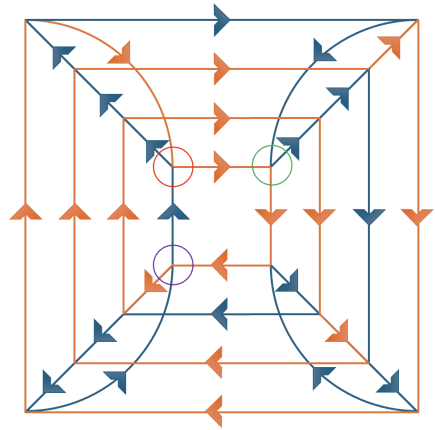}
    \caption{$C_{4} \ \square \ C_{4}$ with Homomorphism and Sample Vertex Types. The sunlet's cycle is in orange and its rays are in blue.}
    \label{fig:s=1}
\end{figure}

Split the vertices into three subsets. The first of these, titled \textbf {Horizontal-Horizontal (H-H)}, denoted in the green circle in Figure \ref{fig:s=1}, are of the points which the sunlet's cycle is drawn from a horizontal edge to another horizontal edges; rays on these points are the edges from $(x-1,y)$ to $(x,y)$. The second subset, \textbf {Horizontal-Vertical (H-V)}, denoted by the purple circle in Figure \ref{fig:s=1}, holds the points which are already connected by horizontal-vertical edge connections and the rays for these points are from $(x-1,y)$ to $(x,y)$. The third subset, \textbf {Vertical-Horizontal (V-H)}, connects the vertical-horizontal edges, denoted by the red circle in Figure \ref{fig:s=1}. These rays lie on the edges from $(x,y-1)$ to $(x,y)$. A visualization of the three subsets is shown in Figure \ref{fig:s=2}.

\begin{figure}[ht]
    \centering
    \includegraphics[scale=.3]{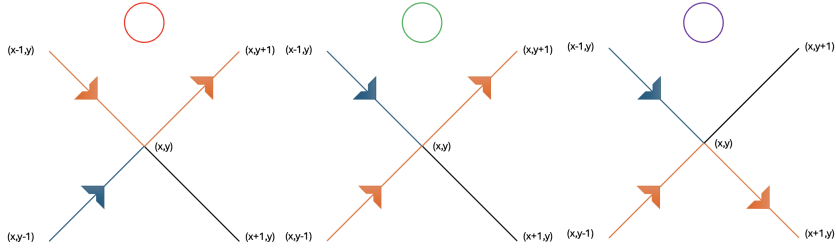}
    \caption{Visualized Point Subsets.}
    \label{fig:s=2}
\end{figure}

For a covering, the $n^{2}$ rays of the sunlet must map to the  $n^{2}$ edges of $C_{n} \ \square \ C_{n}$ not on the Hamiltonian cycle, once and only once through a bijection. Each vertex now only has two remaining edges that are not covered, one with directionality facing in and the other out. 

Following the rules of the three subsets, every vertex's ray follows the remaining edge that points in. Since $n \geq 3$, no vertices have their remaining inward facing edge be the same as their outward facing edge, so the rays into each of the $n^2$ vertices of the sunlet map to the remaining edges  bijectively.

Thus, there is a 2-to-1 covering $\varphi$ of $C_{n} \ \square \ C_{n}$ by $S_{n^2}^{1}$ such that $\varphi(C_{n^2})$ is a spanning cycle and 
$\varphi$ is compatible with FSM and standard orientations.
\end{proof}

\section{Homomorphism from a Set of $n$ Sunlets}

Here we need a bipartite toroidal grid and use the standard orientation.

\begin{thm}
For $n \geq 2$, there is a 2-to-1 covering
$$\varphi:n \cdot S_{4n}^{1} \to C_{2n} \ \square \ C_{2n}$$
such that the restriction $\varphi \,| \, n \cdot Z(S_{4n}^1)$ is an isomorphism onto the disjoint staircases defined below and $\varphi$ is compatible with FSM and standard orientations.
\end{thm}

\begin{figure}[ht]
    \centering
    \includegraphics[scale=.3]{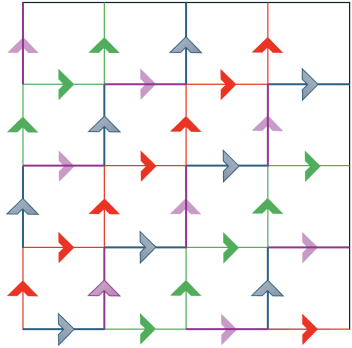}
    \caption{Homomorphism of $2 \cdot S_{8}^{1}$ onto $C_{4} \, \Box \, C_{4}$ on a toroidal grid. One $S_{8}^{1}$ is a red cycle with purple rays,  the other is a green cycle with blue rays; opposite sides identified.} 
    \label{fig:s=3}
\end{figure}

\begin{proof}
A $4n$-cycle in $C_{2n} \ \square \ C_{2n}$ is given as follows:

$$
(00),(10),(11),(21),(22),(32),\ldots, (2n{-}1,2n{-}2), (2n{-}1,2n{-}1),(0,2n{-}1),$$

\noindent where we have written $(ab)$ for $(a,b)$ when it is clear.  We visualize this cycle as a \textbf {(closed) staircase} in the $2n \times 2n$ toroidal grid. There are $n{-}1$ additional closed staircases that are pairwise vertex-disjoint, starting at the points 
$$(2,0), (4,0), \ldots, (2n-2,0).$$

The corresponding function $\varphi: n\, V(C_{4n}) \to V(C_{2n} \ \square \ C_{2n})$ covers each vertex exactly once.
Take the fixed orientation of the grid to be increasing in both coordinate.  There are two forward edges at each vertex in the toroidal grid.  The sunlet cycle mapped to the staircase uses exactly one of the forward directions at a given vertex so we may assign to the ray $uu'$ at cycle vertex $u$ in $S_{4n}^1$ the other forward direction at $\varphi(u)$.  Thus, $\varphi(u')$ is in a different staircase cycle than $\varphi(u)$. At each turn of the staircase, the inward edge is rotated by $\pi/2$, alternating between clockwise and counterclockwise.
 See Fig \ref{fig:s=3}.
\end{proof}

\section{Factorization into Minimum Bipartite Sunlets}

For $n \geq 2$, define the {\bf raster orientation} of  $C_{2n} \ \square \ C_{2n}$ (thought of as a square grid with opposite sides identified) by making all odd-indexed rows left-to-right and all even-indexed rows right-to-left and all odd-indexed columns bottom-to-top and even columns top-to-bottom.  The product of the odd-indexed row and column edges is a set of length-4 dicycles we call the \textbf {odd squares}.

A minimum-vertex bipartite sunlet is a sunlet of form $S_4^1$.

\begin{thm}
For $n \geq 2$, there is a 2-to-1 covering 
$$\varphi: n^2 \cdot S_4^1 \to C_{2n} \,\Box \,C_{2n}$$
such that the restriction $\varphi \,| \,n^2 \cdot Z(S_4^1)$ is an isomorphism onto the odd squares and $\varphi$ is compatible with the FSM and raster orientations.
\end{thm}

\begin{proof}
Let $E_n$ be the subset of the odd-indexed edges of $C_{2n}$.  Then $E_n \,\Box \, E_n$ determines $n^2$ pairwise-disjoint copies of $C_4$.  As above, the complementary graph is a disjoint union of cycles and so can be oriented into dicycles, showing the existence of many choices of where the rays at the vertices of the sunlets are mapped by the homomorphism. But a simple algorithm maps each of the 4 corners of an alternating $C_4$ in the direction of the side terminating at that vertex.  As each vertex of $C_{2n} \,\Box \,C_{2n}$ is in a unique $C_4$, this chosen direction is where to map the ray at that vertex of that copy of $C_4$ as displayed in Figure 5 and Figure 6. At each turn of an odd square, the inward edge is a $\pi/2$ rotation always clockwise.
\end{proof}

\begin{figure}[ht]
    \centering
    \includegraphics[scale=.3]{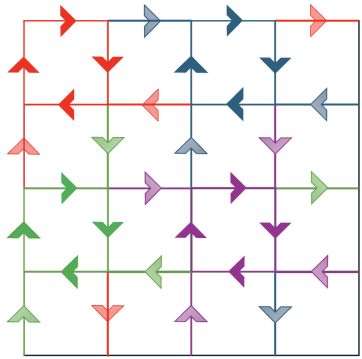}
    \caption{Covering of $4 \cdot S_{4}^{1}$ onto $C_{4}\, \Box\,  C_{4}$ in blue, red, green, and yellow. The cycle edges are denoted by solid arrows and rays by shaded arrows which also indicate directionality.}
    \label{fig:s=5}
\end{figure}

\begin{figure}[ht]
    \centering
    \includegraphics[scale=.3]{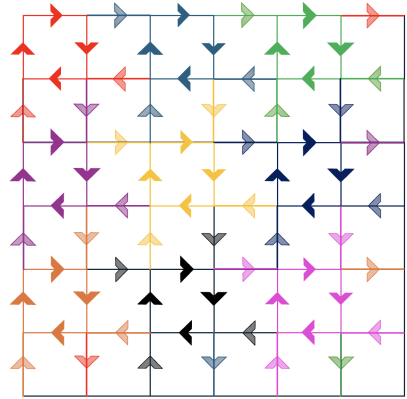}
    \caption{Covering of $9 \cdot S_{4}^{1}$ onto $C_{6}\, \Box\,  C_{6}$ in red, blue, green, purple, yellow, navy, orange, black, and pink. The cycle edges are denoted by solid arrows and rays by shaded arrows which also indicate directionality.}
    \label{fig:s=6}
\end{figure}

\section{Conclusion}
The Cartesian product of two cycles has appeared in computer architectures (e.g., \cite{song}) so it is interesting that it has several rearrangements as an edge-decomposition into sunlets under homomorphisms that are 2-to-1 on vertices.  Sunlets admit exactly two orientations as the digraph of a finite state machine.  Perhaps this could be used to supervise a family of parallel computations that run on the regular toroidal structure.

\end{document}